\documentclass[reqno,a4paper]{amsart}

\usepackage[active]{srcltx}
\usepackage{fancyhdr}

\usepackage[utf8]{inputenc}

\usepackage{times}

\usepackage{latexsym,amsopn,amssymb,amsmath}%{mathabx}
\usepackage{amsthm}
\usepackage{amsfonts,amsbsy,amscd,stmaryrd}
\usepackage{paralist}
\usepackage{hyperref}
\hypersetup{colorlinks=true,linkcolor=Emerald,urlcolor=blue} 

\usepackage[usenames,dvipsnames]{color}

\usepackage{mathrsfs} % for \mathscr{F} for example

\newcommand{\C}{\mathbb{C}}

\newcommand{\R}{\mathbb{R}}

\newcommand{\ca}{\mathcal{A}}

\newcommand{\ce}{\mathcal{E}}

\newcommand{\ch}{\mathcal{H}}

\newcommand{\cR}{\mathcal{R}}

\newcommand{\cv}{\mathcal{V}}

\def\rond{\mathscr}
\newcommand{\ra}{\rond{A}}

\newcommand{\rc}{\rond{C}}

\newcommand{\re}{\rond{E}}

\newcommand{\rk}{\rond{K}}

\def\rmb{\mathrm{b}}

\def\e{\mathrm{e}}

\def\i{\mathrm{i}}

\def\vphi{\varphi}
\def\vkappa{\varkappa}
\def\braket#1#2{\langle{#1}|{#2}\rangle}

\newcommand\slim[1]{\mbox{\rm s-}\!  {\lim_{#1}}}

\def\cbu{\mathcal{C}_{\mathrm b}^{\mathrm u}}

\def\ess{\mathrm{ess}}

\def\pprod{\textstyle \prod}

\long\def\symbolfootnote[#1]#2{\begingroup%
\def\thefootnote{\fnsymbol{footnote}}\footnote[#1]{#2}\endgroup}

\newtheorem{theorem}{Theorem}%[section]

\newtheorem{proposition}{Proposition}

\newtheorem{definition}{Definition}
\newtheorem{remark}{Remark}

\newtheorem{example}{Example}

\renewcommand{\hat}{\widehat}
\renewcommand{\tilde}{\widetilde}

\newcommand{\RR}{\mathbb{R}}
\newcommand{\CC}{\mathbb{C}}
\renewcommand{\SS}{\mathbb{S}}

\def\maB{\mathcal{B}} 
\def\maC{\mathcal{C}} 
            
\def\maE{\mathcal{E}}            
\def\maF{\mathcal{F}}

\newcommand{\oX}{\overline{X}}

\newcommand{\oXZ}{\overline{X/Z}}
\newcommand{\maCX}{\maC(\overline{X})}
\newcommand{\maCXY}{\maC(\overline{X/Y})}

\newcommand\ie{{i.\kern2pt e.\ }}

\title[Essential spectrum]{On the Essential Spectrum of N-Body 
Systems with Asymptotically Homogeneous of Order Zero Interactions}

\begin{document}

\vspace*{-10mm}

\author[V. Georgescu]{Vladimir Georgescu} \address{V. Georgescu,
  D\'epartement de Math\'ematiques, Universit\'e de Cergy-Pontoise,
  95000 Cergy-Pontoise, France}
\email{vladimir.georgescu@math.cnrs.fr}

\author[V. Nistor]{Victor Nistor} \address{V. Nistor, Université de
  Lorraine, UFR MIM, Ile du Saulcy, CS 50128, 57045 METZ Cedex 01, and
  \\ Pennsylvania State University, Math. Dept., University Park, PA
  16802, USA} \email{victor.nistor@univ-lorraine.fr}

\begin{abstract} We overview some of our recent results on the
  essential spectrum of $N$-body Hamiltonians with potentials defined
  by functions that have radial limits at infinity. The results extend
  the HVZ theorem which describes the essential spectrum of usual
  $N$-body Hamiltonians. The proof is based on a careful study of
  algebras generated by potentials and their cross-products. We also
  describe the topology on the spectrum of these algebras, thus
  extending to our setting a result of A. Mageira.  Our techniques
  apply to more general classes of potentials associated to
  translation invariant algebras of bounded uniformly continuous
  functions on a finite dimensional vector space $X$.
\end{abstract}

\maketitle

%%%%%%%%%%%%%%%%%%%%%%%%%%%%%%%%%%
%%%%%%%%%%%%%%%%%%%%%%%%%%%%%%%     
\section*{Introduction}
%%%%%%%%%%%%%%%%%%%%%%%%%%%%%%%%%%
%%%%%%%%%%%%%%%%%%%%%%%%%%%%%%

Let $X$ be a finite dimensional real vector space and, for each linear
subspace $Y$ of $X$, let $V_Y : X/Y \to \RR$ be a Borel function. To
simplify the statements, we assume $V_Y=0$, except for a finite number
of $Y$.  We keep the notation $V_Y$ for the function on $X$ given by
$V_Y\circ\pi_Y$, where $\pi_Y:X\to X/Y$ is the natural map.  In this
paper, we use crossed-products of $C^*$-algebras to study the
essential spectrum of Hamiltonians of the form
\begin{equation}\label{eq.Hamiltonian}
  H \, := \, h(P) + {\textstyle\sum_{Y}} V_Y \,,
\end{equation}
under certain conditions on the potentials $V_Y$. Here $h: X^* \to [0,
\infty [$ is a continuous, proper function and $P$ is the momentum
observable (recall that \emph{proper} means that $
{\lim_{|k|\to\infty}} h(k)=+\infty$).  More precisely, $h(P) =
\maF^{-1} M_h \maF$, where $\maF:L^2(X) \to L^2(X^*)$ is the Fourier
transform and $M_h$ is the operator of multiplication by $h$ (formally
$P=-\i\nabla$).

Operators of this form cover the Hamiltonians that are currently the
most interesting (from a physical point of view) Hamiltonians of
$N$-body systems.  Here are two main examples. In a generalized
version of the non-relativistic case, a scalar product is given on
$X$, so, by taking $h(\xi) = |\xi|^2$, we get $h(P) =\Delta$, the
positive Laplacian. In the simplest relativistic case, $X=(\R^3)^N$
and, writing the momentum $P=(P_1,\dots,P_N)$, we have
$h(P)=\sum_{k=1}^N(P_k^2+m_k^2)^{1/2}$ for some real numbers $m_k$.
We refer to \cite{Derezinski-Gerard} for a thorough introduction to
the subject and study of these systems.

To state our results, we need some notations and a definition.  Let
$\SS_X := (X\smallsetminus \{0\}) / \RR_+$ be \emph{the sphere at
  infinity of $X$}, i.e.\ the set of all half-lines $\hat{a} :=
\RR_+a$ with $\RR_+=\, ]0,\infty[$ and $a \in X\setminus\{0\}$.  A
  function $v:X\to\C$ is said \emph{to have uniform radial limits at
    infinity} if $v(\hat{a}):= {\lim_{r\to\infty}} v(ra)$ exists
  uniformly in $\hat{a}\in\SS_X$. Clearly then $v(\hat{a}) =
  {\lim_{r\to\infty}} v(ra+x) \,,\ \forall x\in X$.
More generally, we are interested in functions $v$ such that $
{\lim_{r\to\infty}} v(ra+x)$ exists for all $x\in X$. The limit may
depend on $x$ and defines a function $\tau_\alpha(v):X\to\C$, where
$\alpha := \hat {a}$.  For example, let us consider $v=V_Y$. Then
$\tau_\alpha(V_Y)(x) = {\lim_{r\to\infty}}
V_Y(r\pi_Y(a)+\pi_Y(x))$. In particular, $\tau_\alpha(V_Y)=V_Y$
whenever $\alpha := \hat{a} \subset Y$ (\ie $a \in Y$). On the other
hand, if $V_Y : X/Y \to \CC$ has uniform radial limits at infinity and
$\hat{a} = \alpha\not\subset Y$, then $\pi_Y(\alpha) := \RR_+\pi_Y(a)
= \widehat{a + Y} \in\SS_{X/Y}$ is well defined and
$\tau_\alpha(V_Y)(x) = V_Y(\pi_Y(\alpha))$ turns out to be a
\emph{constant}.  We denote by $\overline{\cup}_\alpha A_\alpha :=
\overline{\cup_{\alpha} A_\alpha}$ the closure of the union of a given
family of sets. Here is the main result of the paper:

\begin{theorem}\label{th:1}
Assume that each of the functions $V_Y:X/Y\to\R$ is bounded and has
uniform radial limits at infinity. For each $\alpha\in\SS_X$, we
denote
\begin{equation}\label{eq:1}
  \tau_\alpha(H) = h(P) + \sum_{Y} \tau_\alpha(V_Y) = h(P) +
  \sum_{Y\supset\alpha} V_Y + \sum_{Y\not\supset\alpha}
  V_Y(\pi_Y(\alpha)) \, .
\end{equation}
Then $ \sigma(\tau_{\alpha}(H)) = [c_\alpha, \infty)$ for some real
$c_\alpha$ and $\sigma_\ess(H) = \overline{\cup}_{\alpha\in \SS_{X}}
\sigma(\tau_\alpha(H))= [\inf_\alpha c_\alpha, \infty)$.
\end{theorem}

Unbounded potentials are considered in Theorem \ref{th:ess4}.  If all
the radial limits are zero, which is the case of the usual $N$-body
potentials, then the terms corresponding to $\alpha\not\subset Y$ are
dropped in \eqref{eq:1}. Consequently, if $h(P) = \Delta$ is the
non-relativistic kinetic energy, we recover the Hunziker, van Winter,
Zhislin (HVZ) theorem. Descriptions of the essential spectrum of
various classes of self-adjoint operators in terms of limits at
infinity of translates of the operators have already been obtained
before, see for example \cite{Helffer-Mohamed, RRS2004} and \cite{LaS}
(in historical order).  Our approach is based on the ``localization at
infinity'' technique developed in \cite{GI2,GI3} in the context of
crossed-product $C^*$-algebras.

Let us sketch the main idea of this approach. Let $\ca$ be a
translation invariant $C^*$-algebra of bounded uniformly continuous
functions containing the functions which have a limit at infinity and
let $\hat\ca$ be its character space. Note that $\hat\ca$ is a compact
topological space that naturally contains $X$ as an open dense subset
and $\delta(\ca)=\hat\ca\setminus X$ can be thought as a boundary of
$X$at infinity. Then to each self-adjoint operator $H$ affiliated%
\footnote{ A self-adjoint operator $H$ on a Hilbert space $\ch$ is
  said to be \emph{affiliated} to a $C^*$-algebra $\ra$ of operators
  on $\ch$ if one has $(H+\i)^{-1}\in\ra$. Then $\varphi(H)\in\ra$ for
  any $\varphi\in\maC_0(\R)$.  } to the crossed product $\ca\rtimes X$
of $\ca$ by the action of $X$ one may associate a family of
self-adjoint operators $H_\vkappa$ affiliated to $\ca\rtimes X$
indexed by the characters $\vkappa\in \delta(\ca)$. This family
completely describes the image of $H$ (in the sense of affiliated
operators) in the quotient of $\ca\rtimes X$ with respect to the ideal
of compact operators. In particular, the essential spectrum of $H$ is
the closure of the union of the spectra of the operators
$H_\vkappa$. These operators are the \emph{localizations at infinity
  of $H$}, more precisely, \emph{$H_\vkappa$ is the localization of
  $H$ at the point $\vkappa$}.

Once chosen the algebra $\ca$, in order to use these techniques of
this paper concrete situations, one still needs the following: (1) to
have a good description of the character space of the abelian algebra
$\ca$, and (2) to have an efficient criterion for affiliation to the
crossed product $\ca\rtimes X$. In this paper, we explain also how to
achieve these two steps.

\section{Crossed products and localization at infinity}
\label{s:gen}

We begin with some general facts.  Let $\maC_\rmb(X)$ be the algebra
of bounded continuous functions, $\maC_0(X)$ the ideal of functions
vanishing at infinity, $\maC(X^+)=\C+\maC_0(X)$, and $\cbu(X)$ the
subalgebra of bounded uniformly continuous functions on $X$.
If $p\in X^*$ and $q\in X$ then we define unitary operators $S_p$ and
$T_q$ on $L^2(X)$ by
\begin{equation}\label{eq:sit}
	(S_p f)(x) \, = \, \e^{\i p(x)} f(x) \quad \text{ and } 
	\quad (T_q f)(q) \, = \, f(x+q)  \, .
\end{equation}
One sometimes writes $p(x)=\braket{x}{p}$. The symbol $A^{(*)}$ used
below means that the relation should hold for $A$ and $A^*$.

\begin{definition} \label{df:pmlp} We say that a bounded operator $A$
  on $L^2(X)$ {\em has the position-momentum limit property}\ if\,
\begin{equation*}  
	\lim_{p\to0} \, \| [S_p, A] \| \, = \, 0\quad \mbox{\rm\ and } 
	\quad  \lim_{q\to0}\| (T_q-1)A^{(*)} \| \, = \, 0 \, .
\end{equation*}
\end{definition}

It is clear that the set of operators which have the position-momentum
limit property is a $C^*$-algebra. This algebra is equal to the (image
in $\maB(L^2(X))$ of the crossed products $\cbu(X)\rtimes X$ (see
\cite{GI2}). More precisely, if $\ca$ is an arbitrary translation
invariant $C^*$-subalgebra of $\cbu(X)$, then there is a natural
realization of the abstract crossed product $\ca\rtimes X$ as a
$C^*$-algebra of operators on $L^2(X)$ and we do not distinguish the
two algebras. To describe this concrete version of $\ca\rtimes X$, and
for later use, we need some more notations. \label{p:ca}

If $\varphi:X\to\C$ and $\psi:X^*\to\C$ are measurable functions then
$\varphi(Q)$ and $\psi(P)$ are the operators on $L^2(X)$ defined as
follows: $\varphi(Q) := M_\phi$ acts as multiplication by $\varphi$
and $\psi(P)= \maF^{-1}M_\psi\maF$, where $\maF$ is the Fourier
transform $L^2(X) \to L^2(X^*)$ and $M_\psi$ is the operator of
multiplication by $\psi$.  The group $C^*$-algebra of $X$, i.e.\ the
closed subspace of $B(L^2(X))$ generated by convolution with
continuous functions with compact support, will be denoted $C^*(X)$.
The map $\psi\mapsto \psi(P)$ is an isomorphism between $\maC_0(X^*)$
and $C^*(X)$.  With these notations, we have: $\ca\rtimes X$ is the
norm closed linear space of bounded operators on $L^2(X)$ generated by
the products $\varphi(Q)\psi(P)$ with $\varphi\in\ca$ and
$\psi\in\maC_0(X^*)$. In particular, $\ca\rtimes X$ consists of
operators that have the position-momentum limit property.

Our next goal is to recall the definition of localizations at infinity
for operators that have the position-momentum limit property. Let us
fix then an algebra $\ca$ such that $\maC(X^+)\subset\ca$.  As we
explained in the introduction, $\hat\ca$ is a compactification of $X$
and $\delta(\ca)=\hat\ca\setminus X$ is a compact set which can be
characterized as the set of characters $\vkappa$ such that
$\vkappa\vert_{\maC_0(X)}=0$.  The largest $\ca$ that we may consider
is $\cbu(X)$ and we denote $\gamma(X) := \widehat{\cbu(X)}$ its
character space and $\delta(X)=\gamma(X)\setminus X$.  For a general
$\ca$ now, the map $\vkappa\mapsto\vkappa\vert_{\ca}$ is a continuous
surjective map $\delta(X)\to\delta(\ca)$, and hence $\delta(\ca)$ is a
quotient of the compact space $\delta(X)$.  If $q\in X$ and $\varphi$
is a function on $X$ then $T_q\vphi$ is its translation by $q$,
defined as in \eqref{eq:sit}. We extend now this definition of $T_q$
by replacing $q \in X$ with a character $\vkappa \in \hat\ca$:
\begin{equation}\label{eq:trinf}
  (T_\vkappa\varphi)(x)= \vkappa(T_x\varphi) \, ,  \quad\text{for any }\
  \varphi\in\ca \,, \ \ \vkappa\in\hat\ca\ \text{\ \ and, }\ x\in X . 
\end{equation}
It is clear that $T_\vkappa\varphi\in\cbu(X)$ and that its definition
coincides with the previous one if $\vkappa=q\in X$.  Moreover, we
also get ``translations at infinity'' of $\varphi\in\ca$ by elements
$\vkappa\in\delta(\ca)$; note however that such a translation does not
belong to $\ca$ in general. Also, the function $\vkappa\mapsto
T_\vkappa\varphi\in \cbu(X)$ defined on $\hat\ca$ is continuous if
$\cbu(X)$ is equipped with the topology of local uniform convergence,
hence $T_\vkappa\varphi = {\lim_{q\to\vkappa}}\, T_q\vphi$ in this
topology for any $\vkappa\in\delta(\ca)$.  If $A$ is an operator on
$L^2(X)$, let $\tau_q(A)=T_qAT_q^*$ be its translation by $q\in
X$. Clearly $\tau_q(\varphi(Q))=(T_q\varphi)(Q)$.  If $A\in\ca\rtimes
X$, then we may also consider ``translations at infinity'' by elements
of the boundary $\delta(\ca)$ of $X$ in $\hat\ca$ and we get a useful
characterization of the compact operators.  The following facts are
mainly consequences of \cite[Theorem 1.15]{GI3}:
\begin{itemize}

\item If $\vkappa\in\hat\ca$, then there is a unique morphism
  $\tau_\vkappa:\ca\rtimes X\to\cbu(X)\rtimes X$ such that
\begin{equation*}
   \tau_\vkappa(\varphi(Q) \psi(P) ) \, = \,
   (T_\vkappa\varphi)(Q)\psi(P) \,, \quad\text{for all }
   \varphi\in\cbu(X), \, \psi\in\maC_0(X) \, .
\end{equation*}

\item If $A\in\ca\rtimes X$ then $\vkappa\mapsto \tau_\vkappa(A)$ is a
  strongly continuous map $\hat\ca\to \maB(L^2(X))$.

\item $\cap_{\vkappa\in\delta(\ca)}\ker\tau_\vkappa=\maC_0(X)\rtimes
  X\equiv\rk(X)=$ ideal of compact operators on $L^2(X)$.

\item The map $\tau(A)=(\tau_\vkappa(A))_{\vkappa\in\delta(\ca)}$
  induces an injective morphism
\begin{equation}\label{eq:crk}
  	\ca\rtimes X/\rk(X) \hookrightarrow
        \pprod_{\vkappa\in\delta(\ca)}\rc(X) .
\end{equation}

\item If $A\in\ca\rtimes X$ is a normal operator then $\sigma_\ess(A)
  = \overline{\cup}_{\vkappa\in\delta(\ca)} \sigma (\tau_\vkappa(A))$.
\item If $H$ is a self-adjoint operator on $L^2(X)$ affiliated to
  $\ca$ then for each $\vkappa\in\delta(\ca)$ the limit
  $\tau_\vkappa(H):=\slim{q\to\vkappa}T_qHT_q^*$ exists and
  $\sigma_\ess(H) = \overline{\cup}_{\vkappa\in\delta(\ca)} \sigma
  (\tau_\vkappa(H))$.
\end{itemize}

To be precise, the last strong limit means: $\tau_\vkappa(H)$ is a
self-adjoint operator (not necessarily densely defined) on $L^2(X)$
and $\slim{q\to\vkappa}\theta(T_qHT_q^*)=\theta(\tau_\vkappa(H))$ for
all $\theta\in\maC_0(\R)$.  It is clear that in the last three
statements above one may replace $\delta(\ca)$ by a subset $\pi$ if
for each $A\in\ca\rtimes X$ we have: $\tau_\vkappa(A)
=0\,\forall\vkappa\in\pi\Rightarrow
\tau_\vkappa(A)=0\,\forall\vkappa\in\delta(\ca)$.  In the case of
groupoid (pseudo)differential algebras (that is, when $\hat \ca$ is a
manifold with corners), the morphisms $\tau_\vkappa$ can be defined
using restrictions to fibers, as in \cite{LMN1}, and the last three
statements above remain valid.

\section{The Two-Body Case}\label{s:2b}

As a warm up and in order to introduce some general notation, we treat
first the two-body case, where complete results may be obtained by
direct arguments.  The algebra of interactions in the standard
two-body case is $\maC(X^+)$, and hence the Hamiltonian algebra is
\begin{equation}\label{eq:s2b}
	\maC(X^+)\rtimes X = \C\rtimes X + \maC_0(X)\rtimes X
        =C^*(X)+\rk(X)
\end{equation}
where the sums are direct. Thus $\maC(X^+) \rtimes X/\rk(X)=C^*(X)$
which finishes the theory.  Another elementary case, which has been
considered as an example in \cite{GI2}, is $X=\R$ with $\maC(\R^+)$
replaced by the algebra $\maC(\overline{\R})$ of continuous functions
which have limits (distinct in general) at $\pm\infty$. Then there is
no natural direct sum decomposition of $\maC(\overline{\R})\rtimes \R$
as in \eqref{eq:s2b}, but one has by standard arguments
\begin{equation}\label{eq:h2b}
  	\maC(\overline{\R})\rtimes\R / \rk(\R) \simeq C^*(\R) \oplus
        C^*(\R).
\end{equation}
Our purpose in this section is to extend \eqref{eq:h2b} to arbitrary
$X$ and the main result is an explicit description of the crossed
product algebra $\maC(\overline{X})\rtimes X$.

Let $\maCX$ be the closure in $\maC_\rmb(X)$ of the subalgebra of
functions homogeneous of degree zero outside a compact set.  Then $
\maCX = \{ u\in \maC(X) \mid {\lim_{\lambda\to+\infty}} u(\lambda a)
\ \text{ exists uniformly in }\ \hat a\in\SS_X \} ,$ where, we recall,
$\hat a := \RR_+ a$ and $\SS_X := (X \smallsetminus \{0\})/\RR_+$, so
$\hat{a} \in \SS_X$.  As a set, the character space of $\maC(\oX)$ can
be identified with the disjoint union $\oX = X \cup \SS_{X}$. The
topology induced by the character space on $X$ is the usual one and
the intersections with $X$ of the neighbourhoods of some
$\alpha\in\SS_X$ are the sets which contain a truncated cone $C$ such
that there is $a\in\alpha$ such $\lambda a \in C$ if
$\lambda\geq1$. The set of such subsets is a filter $\tilde\alpha$ on
$X$ and, if $Y$ is a Hausdorff space and $u:X\to Y$, then
$\lim_{\tilde\alpha} u=y$ means that $u^{-1}(V)\in\tilde\alpha$ for
any neighborhood $V$ of $u$.  We shall write $ {\lim_{x\to\alpha}}
u(x)$ instead of $\lim_{\tilde\alpha} u$.  We have that $\maCX$ is a
translation invariant $C^*$-subalgebra of $\cbu(X)$ and so the crossed
product $\maCX\rtimes X$ is well defined. We have the following
explicit description of this algebra.

\begin{proposition}\label{pr:sphalg}
The algebra $\maCX\rtimes X$ acting on $L^2(X)$ consists of bounded
operators $A$ that have the position-momentum limit property and are
such that the limit $\tau_\alpha(A) = \slim{a\to\alpha} T_a A T_a^*$
exists for each $\alpha = \hat a \in\SS_X$.  If $A\in\maCX\rtimes X$
and $\alpha\in\SS_X$, then $\tau_\alpha(A) \in C^*(X)$ and
$\tau(A):\alpha\mapsto\tau_\alpha(A)$ is norm continuous. The map
$\tau: \maCX\rtimes X \to C(\SS_X)\otimes C^*(X)$ is a surjective
morphism whose kernel is the set of compact operators on $L^2(X)$,
which gives
\begin{equation}\label{eq:quot}
	\maCX\rtimes X /\rk(X) \cong C(\SS_X)\otimes C^*(X) .
\end{equation}
If $H$ is a self-adjoint operator affiliated to $\maCX\rtimes X$ then
$\tau_\alpha(H)=\slim{a\to\alpha} T_aHT_a^* $ exists for all
$\alpha\in\SS_X$ and $\sigma_\ess(H)=\cup_\alpha
\sigma(\tau_\alpha(H))$.
\end{proposition}

Equations \eqref{eq:h2b} and \eqref{eq:quot} follow also from standard
properties of crossed-product algebras by amenable groups.

In the next two example $H=h(P)+V$ with $h:X^*\to[0,\infty[$
continuous and proper. We denote by $|\cdot|$ a fixed norm on $X^*$
and by $\ch^s$ we denote the usual Sobolev spaces on $X$ ($s\in\R$).

\begin{example}\label{ex:2.1}{\rm 
Let $V$ be a bounded symmetric operator satisfying: (1) $
\lim_{p\to0}\| [S_p, V] \|=0$ and (2) the limit
$\tau_\alpha(V)=\slim{a\to\alpha} T_a V T_a^*$ exists for each
$\alpha\in\SS_X$.  Then $H$ is affiliated to $\maCX\rtimes X$ and
$\tau_\alpha(H)=h(P)+\tau_\alpha(V)$. Moreover, if $V$ is a function
then $\tau_\alpha(V) $ is a number, but in general we have
$\tau_\alpha(V)=v_\alpha(P)$ for some function
$v_\alpha\in\cbu(X^*)$.}
\end{example}

\begin{example}\label{ex:2.2}
{\rm Assume that $h$ is locally Lipschitz and that there exist $c,s>0$
  such that, for all $p$ with $|p|>1$,
\begin{equation}\label{eq:ess2}
  |\nabla h (p)|\leq c \big( 1+h(p) \big) \quad\text{and}\quad
  c^{-1}|p|^{s} \leq \big( 1+h(p) \big)^{1/2} \leq c |p|^{s} .
\end{equation}
Let $V:\ch^s\to\ch^{-s}$ such that $\pm V\leq \mu h(P) +\nu$ for some
numbers $\mu,\nu$ with $\mu<1$ and satisfying the next two conditions:
(1) $ {\lim_{p\to0}} \| [S_p, V] \|_{\ch^s\to\ch^{-s}}=0$, (2)
$\forall\alpha\in\SS_X$ the limit $\tau_\alpha(V)=\slim{a\to\alpha}
T_a V T_a^*$ exists strongly in $\maB(\ch^s,\ch^{-s})$.  Then $h(P)+V$
and $h(P)+\tau_\alpha(V)$ are symmetric operators $\ch^s\to\ch^{-s}$
which induce self-adjoint operators $H$ and $\tau_\alpha(H)$ in
$L^2(X)$ affiliated to $\maCX\rtimes X$ and
$\sigma_\ess(H)=\cup_\alpha \sigma(\tau_\alpha(H))$.  }
\end{example}

\section{N-body case}\label{sec.Nbody}

We first indicate a general way of constructing $N$-body
Hamiltonians. For each linear subspace $Y\subset X$ let
$\ca(X/Y)\subset\cbu(X/Y)$ be a translation invariant $C^*$-subalgebra
containing $\maC_0(X/Y)$ with $\ca(X/X)=\ca(0)=\C$. We embed
$\ca(X/Y)\subset\cbu(X)$ as usual by identifying $v$ with
$v\circ\pi_Y$. Then the $C^*$-algebra $\ca$ generated by these
algebras is a translation invariant $C^*$-subalgebra of $\cbu(X)$
containing $\maC(X^+)$ and thus we may consider the crossed product
$\ca\rtimes X$ which is equal to the $C^*$-algebra generated by the
crossed products $\ca(X/Y)\rtimes X$.  The operators affiliated to
$\ca\rtimes X$ are $N$-body Hamiltonians.

The standard $N$-body algebra corresponds to the minimal choice
$\ca(X/Y)=\maC_0(X/Y)$ and has remarkable properties which makes its
study relatively easy (it is graded by the lattice of subspaces of
$X$). Our purpose in this paper is to study what could arguably be
considered to be the simplest extension of the classical $N$-body
setting: namely, it is obtained by choosing $\ca(X/Y)=\maCXY$ for all
$Y$. The next more general case would correspond to the choice
$\ca(X/Y)=\cv(X/Y)$ which is the algebra of functions with vanishing
oscillating at infinity (the closure in $\sup$ norm of the set of
bounded functions of class $C^1$ with derivatives tending to zero at
infinity).

\begin{definition}\label{df:main}
Let $\ce(X)$ be the $C^*$-subalgebra of $\cbu(X)$ generated by $\cup_Y
\maCXY$.
\end{definition}

Clearly $\ce(X)$ is a translation invariant $C^*$-subalgebra of
$\cbu(X)$ containing $\maC(X^+) := \maC_0(X) + \CC$.  If $Y$ is a
linear subspace of $X$ then the $C^*$-algebra
$\ce(X/Y)\subset\cbu(X/Y)$ is well defined and naturally embedded in
$\ce(X)$: it is the $C^*$-algebra generated by $\cup_{Z\supset Y}
C(\oXZ)$. We have
\begin{equation}\label{eq:cexyz}
  \C=\ce(0)=\ce(X/X)\subset\ce(X/Y)\subset\ce(X/Z)\subset\ce(X).
\end{equation}
If $\alpha\in\SS_X$, we shall denote by abuse of notation $X/\alpha$
be the quotient $X/[\alpha]$ of $X$ by the subspace $[\alpha] := \RR
\alpha$ generated by $\alpha$ and let us set
$\pi_\alpha=\pi_{[\alpha]}$. It is clear that $\tau_\alpha(u)(x) =
{\lim_{r\to+\infty}} u(ra+x)$ exists $\forall u\in\ce(X)$ and that the
resulting function $\tau_\alpha(u)$ belongs to $\ce(X)$. The map
$\tau_\alpha$ is an endomorphism of $\ce(X)$ and a linear projection
of $\ce(X)\,$ onto the $C^*$-subalgebra $\ce(X/\alpha)$.

If $\alpha\in\SS_X$ and $\beta\in\SS_{X/\alpha}$, then $\beta$
generates a one dimensional linear subspace $[\beta]:= \RR \beta
\subset X/\alpha$, as above, and hence $\pi_\alpha^{-1}([\beta])$ is a
two dimensional subspace of $X$ that we shall denote
$[\alpha,\beta]$. We shall identify $(X/\alpha)/\beta$ with
$X/[\alpha,\beta]$. Then we have two idempotent morphisms
$\tau_\alpha:\ce(X)\to\ce(X/\alpha)$ and
$\tau_\beta:\ce(X/\alpha)\to\ce(X/[\alpha,\beta])$. Thus
$\tau_\beta\tau_\alpha:\ce(X)\to\ce(X/[\alpha,\beta])$ is an
idempotent morphism.  This construction extends in an obvious way to
families $\overrightarrow\alpha=(\alpha_1,\dots,\alpha_n)$ with
$n\leq\dim X$ and $\alpha_1\in\SS_X, \alpha_2\in\SS_{X/\alpha_1},
\alpha_3\in\SS_{X/[\alpha_1,\alpha_2]}$, \dots (we allow $n=0$ and
denote $A$ the set of all such families). The endomorphism
$\tau_{\overrightarrow\alpha}$ of $\ce(X)$ is defined by induction:
$\tau_{\overrightarrow\alpha}=\tau_{\alpha_n}\dots \tau_{\alpha_1}$.
We also define $[\alpha_1,\alpha_2,\ldots,\alpha_n]$ by induction, so
this is an $n$-dimensional subspace of $X$ associated to $(\alpha_1,
\alpha_2, \ldots, \alpha_n)$ and we denote $X/\overrightarrow\alpha$
the quotient of $X$ with respect to it.  Thus
$\tau_{\overrightarrow\alpha}$ is an endomorphism of $\ce(X)$ and a
projection of $\ce(X)$ onto the subalgebra
$\ce(X/\overrightarrow\alpha)$.  The next description of the spectrum
of $\ce(X)$ extends some results of A. Mageira \cite{Mageira}.

\begin{theorem}\label{th:char}
If $\overrightarrow\alpha\in A$ and $a\in X/\overrightarrow\alpha$,
then $\vkappa(u)=(\tau_{\overrightarrow\alpha}u)(a)$ defines a
character of $\ce(X)$. Conversely, each character of $\ce(X)$ is of
this form.
\end{theorem}

\begin{remark}\label{re:open}{\rm A natural abelian $C^*$-algebra in
the present context is the set $\cR(X)$ of all bounded uniformly
continuous functions $v:X\to\C$ such that $\lim_{r\to\infty} v(ra+x)$
exists, for each $a\in X$, locally uniformly in $x\in X$. This algebra
is larger than $\ce(X)$. It would be interesting to find an explicit
description of its spectrum.  }
\end{remark}

We may now state our main results: they are consequences of the
general theory developed so far (Section \ref{s:gen}, Proposition
\ref{pr:sphalg}, Theorem \ref{th:char}, and Example
\ref{ex:2.2}).

\begin{theorem}\label{th:ess} 
Let $H$ be a self-adjoint operator on $L^2(X)$ affiliated to
$\ce(X)\rtimes X$.  Then $\slim{r\to+\infty} T_{ra} H
T_{ra}^*=:\tau_{\hat{a}}(H)$ exists $\forall a\in X\setminus\{0\}$ and
$\sigma_\ess(H) \, = \, \overline{\cup}_{\alpha\in\SS_X}
\sigma(\tau_\alpha(H)).$
\end{theorem}

\begin{theorem}\label{th:ess4}
Let $h$ be as in Example \ref{ex:2.2} and $V=\sum V_Y$ with
$V_Y:\ch^s\to\ch^{-s}$ symmetric operators such that $V_Y=0$ but for a
finite number of $Y$ and satisfying:
\begin{compactenum}
\item $\exists\,\mu_Y,\nu_Y\geq0$ with $\sum_Y\mu_Y<1$ such that $\pm
  V_Y\leq \mu_Y h(P) +\nu_Y$,
\item $ \lim_{p\to0}\| [S_p, V_Y]   \|_{\ch^s\to\ch^{-s}}=0$,
\item $[T_y,V_Y]=0$ for all $y\in Y$,
\item $\tau_\alpha(V_Y):=\slim{a\to\alpha} T_a V_Y T_a^*$
  exists in $B(\ch^s,\ch^{-s})$ for all $\alpha\in\SS_{X/Y}$.
\end{compactenum}
Then the maps $\ch^s\to\ch^{-s}$ given by $h(P)+V$ and
$h(P)+\sum_Y\tau_\alpha(V_Y)$ induce self-adjoint operators $H$ and
$\tau_\alpha(H)$ in $L^2(X)$ affiliated to $\re(X)$ and
$\sigma_\ess(H) \, = \, \overline{\cup}_{\alpha\in\SS_X}
\sigma(\tau_\alpha(H))$.
\end{theorem} 

Using also some techniques from \cite{DG1}, we obtain:

\begin{example}\label{ex:4}{\rm 
Theorem \ref{th:ess4} covers operators $H=\sum_{|k|,|\ell|\leq s} P^k
a_{k\ell}P^\ell$ which are uniformly elliptic with the $a_{k\ell}$
finite sums of functions of the form $v_Y\circ\pi_Y$ with
$v_Y:X/Y\to\R$ bounded measurable and such that $ {\lim_{z\to\alpha}}
v_Y(z)$ exists uniformly in $\alpha\in\SS_{X/Y}$. The fact that we
allow $a_{k\ell}$ to be only bounded measurable in the principal part
of the operator (\ie the terms with $|k|=|\ell|=s$) is not trivial.  }
\end{example}

In addition to the above mentioned results, we also use general
results on cross-product $C^*$-algebras, their ideals, and their
representations \cite{Dixmier, Skandalis, Renault}.  The maximal ideal
spectrum of the algebra $\maE(X)$ is of independent interest and can
be used to study the regularity properties of the eigenvalues of the
$N$-body Hamiltonian \cite{ACN}. Its relation to the constructions of
Vasy in \cite{VasyReg} will be studied elsewhere. 

{\bf Acknowledgements.} We thank Bernd Ammann for several
  useful discussions.

\end{document}